\documentclass[12pt,fleqn,leqno]{article}
\usepackage{amsmath,amsfonts}

\let\section\undefined
\let\section=\subsection

\def\thmcounterend{.}
\renewenvironment{@begintheorem}[2]{\sl \trivlist
   \item[\hskip \labelsep{\bf #1\ #2\thmcounterend}]}

\newenvironment{proof*}[1]{\trivlist
   \item[\hskip \labelsep\textbf{#1.} ]}
{{\hfill\vbox{\hrule\hbox{%
\vrule height1.3ex\hskip1ex\vrule}\hrule
}}\endtrivlist}

\def\amsclass{\vskip 1\baselineskip
  \typeout{AMS Classification}
    {\noindent\textbf{2000 Mathematical
        Subject Classification: }}}

\def\keywords{\vskip 1\baselineskip
  \typeout{Key words}
    {\noindent\textbf{Key words: }}}

\newcounter{rmnum}

\newtheorem{theorem}{Theorem}[subsection]
\newtheorem{lemma}[theorem]{Lemma}
\newtheorem{corollary}[theorem]{Corollary}
\newtheorem{proposition}[theorem]{Proposition}

\newtheorem{claim}{Claim}

\newcommand{\R}{\mathbb{R}}
\newcommand{\N}{\mathbb{N}}
\newcommand{\ball}{\mathbb{B}}
\newcommand{\sphere}{\mathbb{S}}
\newcommand{\uin}{\mathbb{I}}
\newcommand{\bfd}{\boldsymbol{d}}

\newcommand{\diam}{\hbox{\rm diam}}

\newcommand{\conv}{\hbox{\rm conv}}

\newcommand{\cone}{\hbox{\rm Cone}}
\newcommand{\graph}{\hbox{\rm Graph}}

\newcommand{\sel}{\mathcal{S}\hspace{-1pt}e\ell}

\begin{document}

\title{Dense families of selections and finite-dimensional spaces}

\author{Valentin Gutev and Vesko Valov}

\date{}

\maketitle
\begin{abstract}
  A characterization of $n$-dimensional spaces via continuous
  selections avoiding $Z_n$-sets is given, and a selection theorem for
  strongly countable-dimensional spaces is established.  We apply
  these results to prove a generalized Ostrand's theorem, and to
  obtain a new alternative proof of the Hurewicz formula. It is also
  shown that our selection theorem yields an easy proof of a Michael's
  result.
  \amsclass 54C60, 54C65, 55M10.
  \keywords{Continuous selection, finite-dimensional space,
  strongly countable-dimensional space, $Z_n$-set.}
\end{abstract}

\renewenvironment{@opargbegintheorem}[3]{\sl \trivlist
\item[\hskip \labelsep{\bf #1\ #2\ (#3)\thmcounterend}]}

\section{Introduction}

One of the best characterizations of the covering dimension
is given by extensions of maps in spheres of Euclidean
spaces. Namely, a normal space $X$ has a covering dimension
$\dim(X)\le n$ iff for every closed $A\subset X$, every
continuous map $g:A\to \sphere^{n}$ can be extended to a
continuous map $f:X\to \sphere^{n}$. Here, $\sphere^{n}$
denotes the \emph{$n$-sphere}.\medskip

The parametric version of this fact is provided by the following
Michael's selection theorem in \cite{michael:56b} (see, also,
\cite{repovs-semenov:98}).

\begin{theorem}[\cite{michael:56b}]\label{th:motiv-1}
  Let $X$ be a paracompact space, with $\dim(X)\le n$, $Y$ be a
  completely metrizable space, and let $\varphi :X\rightarrow
  \mathcal{F}(Y)$ be an l.s.c.\ mapping such that $\{\varphi(x):x\in
  X\}$ is equi-$LC^{n-1}$ in $Y$ and each $\varphi(x)$, $x\in X$, is
  $C^{n-1}$. Then, $\varphi$ has a continuous selection.
\end{theorem}

Here, $\mathcal{F}(Y)=\{S\in 2^{Y}: S\ \hbox{is closed}\}$, where
$2^Y$ is the family of all non-empty subsets of $Y$.  A set-valued
mapping $\varphi :X\rightarrow 2^Y$ is \emph{lower semi-continuous},
or l.s.c., if $\varphi ^{-1}(U)=\{x\in X: \varphi (x)\cap
U\neq\emptyset\}$ is open in $X$ for every open $U\subset Y$, and a
map $f:X\to Y$ is a \emph{selection} for $\varphi:X\to 2^Y$ if
$f(x)\in\varphi(x)$ for every $x\in X$. \medskip

Let $m\geq -1$. A family ${\cal S}$ of subsets of a space $Y$ is
\emph{equi-$LC^m$ in $Y$} \cite{michael:56b} if, for every $y\in Y$
and neighbourhood $U$ of $y$, there exists neighbourhood $V$ of $y$
such that, for every $S\in {\cal S}$, every continuous image of
$\sphere^{k}$ ($k\leq m$) in $V\cap S$ is contractible in $U\cap S$.
A space $S$ is called $C^m$ if every continuous image of $\sphere^{k}$
($k\leq m$) in $S$ is contractible in $S$.\medskip

In the present paper we deal with a similar problem but now concerning
another property of the covering dimension.  Towards this end, for
spaces $X$ and $Y$, we use $C(X,Y)$ to denote the set of all
continuous map from $X$ to $Y$. Let $Y$ be a metrizable space, and let
us recall that a closed set $F\subset Y$ (possibly empty) is called a
\emph{$Z_n$-set} in $Y$, $n<\omega$, if the set
$C(\mathbb{B}^n,Y\backslash F)$ is dense in $C(\ball^n,Y)$ with
respect to the uniform topology generated by a metric on $Y$, see
\cite{chigogidze:96} and \cite{torunczyk:78}. Here, $\ball^n$ denotes
the $n$-dimensional closed ball, and $\omega$ is the first infinite
cardinal.  Finally, we shall say that $F$ is a \emph{$\sigma Z_n$-set}
in $Y$ if $F$ is a countable union of $Z_n$-subsets of $Y$.  \medskip

The $\sigma Z_{n}$-sets are ``thin''-subsets with respect to maps from
$n$-dimensional spaces. Namely, a normal space $X$ has a covering
dimension $\dim(X)\le n$ iff for every $\sigma Z_{n}$-set $F$ in
$\ball^{n+1}$, the set $C(X,\ball^{n+1}\setminus F)$ is dense in
$C(X,\ball^{n+1})$ with respect to the uniform topology, see
\cite{epol:88b}.  Our main goal is now to provide a parametric version
of the above characterization of $\dim(X)\le n$. To state it, we need
a bit more terminology about set-valued mappings.  Let
$\mathcal{P}(Y)=2^{Y}\cup \{\emptyset\}$.  To every mapping $\phi:X\to
\mathcal{P}(Y)$ we associate the following subsets of $C(X,Y)$:
\[
\sel(\phi)=\{f\in C(X,Y): f(x)\in\phi (x)\,\ \hbox{for every}\ x\in X\},
\]
and
\[
\mathcal{M}(\phi)=\{f\in C(X,Y): f(x)\notin \phi(x)\,\ \hbox{for
  every}\ x\in X\}.
\]
Also, let us recall that $\phi:X\to \mathcal{P}(Y)$ has a \emph{closed
  graph} (respectively, an \textit{$F_{\sigma}$-graph}) if
$\graph(\phi)=\{(x,y)\in X\times Y: y\in\phi(x)\}$ is closed
(respectively, $F_{\sigma}$) in $X\times Y$.\medskip

Finally, let us recall that, for a metric space $(Y,\rho)$, the
\emph{fine topology} on $C(X,Y)$ is the topology in which the family
of all sets
\[
V(f,\alpha)=\{g\in C(X,Y):\rho(g(x),f(x))<\alpha (x)\ \hbox{for each
  $x\in X$}\},
\]
is a local base at $f$ (see \cite{munkres:75}), where $\alpha$ runs on
the positive continuous functions on $X$. For any space $X$ the fine
topology is finer than the uniform one (generated by $\rho$), and it
does not depend on the metric of $Y$ provided $X$ is normal and
countably paracompact \cite{dimaio-hola-holy-mccoy:98} (for a
paracompact $X$, see \cite{krikorian:69}).\medskip

The following theorem will be proven in this paper.

\begin{theorem}\label{th:fdver}
  Let $X$ be a paracompact space, with $\dim(X)\le n$, $Y$ be a
  completely metrizable space, $\varphi :X\rightarrow \mathcal{F}(Y)$
  be an l.s.c.\ mapping such that $\{\varphi(x):x\in X\}$ is
  equi-$LC^{n-1}$ in $Y$ and each $\varphi(x)$, $x\in X$, is
  $C^{n-1}$, and let $\psi:X\to \mathcal{P}(Y)$ be a mapping with an
  $F_{\sigma}$-graph such that $\psi(x)\cap\varphi (x)$ is a $\sigma
  Z_{n}$-set in $\varphi (x)$ for every $x\in X$. Then, the set
  $\mathcal{M}(\psi)\cap \sel(\varphi)$ is a dense $G_{\delta}$-subset
  of $\sel(\varphi)$ with respect to the fine topology on
  $\sel(\varphi)$.
\end{theorem}

A few words about the paper. Sections 2 and 3 contain the preparation
for the proof of Theorem \ref{th:fdver} which will be finally
accomplished in Section 4. The rest of the paper is devoted to some
possible applications.  For instance, in Section 5 we obtain a
generalization of the classical Ostrand's theorem \cite{ostrand:71},
see Theorem \ref{th:ostrand}. In the same Section 5, we apply Theorem
\ref{th:fdver} to establish some properties of strongly
countable-dimensional spaces in the terms of the Baire category which
are similar to some results of E.\ Pol in \cite{epol:88b} and
\cite{epol:88a}.  We also provide an alternative proof of the Hurewicz
formula that $\dim (X)\leq\dim f+\dim (Y)$, where $f\colon X\to Y$ is
a closed continuous map and $\dim f=\sup\{\dim (f^{-1}(y)):y\in Y\}$,
see Theorem \ref{th;hurewicz}.

\section{Selections avoiding closed sets}

Let $(Y,d)$ be a metric space. If $S\subset Y$ and $\varepsilon>0$,
then we use $B_{\varepsilon}^{d}(S)$ to denote the
\emph{$\varepsilon$-neighbourhood} of $S$ in $(Y,d)$, i.e.\
$B_{\varepsilon}^{d}(S)=\{y\in Y: d(y,S)<\varepsilon\}$. Also, we use
$\diam_{d}(S)$ to denote the diameter of $S$ with respect to $d$,
i.e.\ the \emph{$d$-diameter} of $S$ in $Y$. Now, let us recall that a
family $\mathcal{S}$ of subsets of $Y$ is \emph{$d$-uniformly
  equi-$LC^{n}$} (for some $n\ge -1$) if for every $\varepsilon>0$
there exists $\delta(\varepsilon)>0$ such that, for every
$S\in\mathcal{S}$, every continuous image of a $k$-sphere ($k\le n$)
in $S$ of $d$-diameter$\ < \delta(\varepsilon)$ is contractible over a
subset of $S$ of $d$-diameter$\ <\varepsilon$, \cite{michael:56b}.
\medskip

In what follows, for a metrizable space $Y$, let $\mathcal{D}(Y)$ be
the set of all metrics on $Y$ compatible with the topology of $Y$.  If
$d,\rho\in \mathcal{D}(Y)$, then we will use $\rho\ge d$ to denote
that $\rho(y,z)\ge d(y,z)$ for every $y,z\in Y$.  \medskip

In this section, we provide a slight generalization of \cite[Theorem
4.1]{michael:56b} constructing controlled selections avoiding a given
closed set (Theorem \ref{th:avoid-open}). To prepare for this, we need
the following improvement in \cite[Lemma 11.2]{michael:56b}.

\begin{lemma}\label{lm:functions}
  Let $Y$ be a completely metrizable space, $d\in \mathcal{D}(Y)$, and
  let $\mathcal{S}\subset 2^{Y}$ be $d$-uniformly equi-$LC^{n}$. Also,
  let $\rho\in \mathcal{D}(Y)$ be a complete metric on $Y$ such that
  $\rho\ge d$ and $\mathcal{S}'=\mathcal{S}\cup \{Y\}$ is
  $\rho$-uniformly equi-$LC^{n}$.  Then, to every $\varepsilon>0$ and
  $\mu>0$ there corresponds an $\eta(\varepsilon)>0$ and
  $\lambda(\varepsilon,\mu)>0$ with the following property\/$:$ If
  $k\le n$ and $S\in \mathcal{S}$, then every continuous
  $p:\sphere^{k}\to B_{\lambda(\varepsilon,\mu)}^{\rho}(S)$, with
  $\diam_{d}(p(\sphere^{k}))<\eta(\varepsilon)$, is homotopic to a
  constant map over a subset of $B_{\mu}^{\rho}(S)$ of $d$-diameter$\
  <\varepsilon$.
\end{lemma}

\begin{proof}
  The proof almost repeats that of \cite[Lemma 11.2]{michael:56b}.
  Namely, let $\gamma(\varepsilon)\le \varepsilon$ (respectively,
  $\kappa(\mu)$) be as in \cite[Corollary 4.2]{michael:56b}
  (respectively, as in \cite[Lemma 11.1]{michael:56b}) applied to
  $\mathcal{S}'$ as a family being $\rho$-uniformly equi-$LC^{n}$.
  Next, let $\delta(\varepsilon)\le \varepsilon$ be as in the
  definition of $d$-uniformly equi-$LC^{n}$ of $\mathcal{S}$. Finally,
  following \cite[Lemma 11.2]{michael:56b}, define
  $\eta(\varepsilon)=\delta(\varepsilon/2)/2$ and
  $\lambda(\varepsilon,\mu)=\kappa(\gamma(\xi))$, where
  $\xi=\min\{\delta(\varepsilon/2)/4, \mu\}$. This works in our
  present situation. Indeed, let $k\le n$, $S\in \mathcal{S}$, and let
  $p:\sphere^{k}\to B_{\lambda(\varepsilon,\mu)}^{\rho}(S)$ be
  continuous such that $\diam_{d}(p(\sphere^{k}))<\eta(\varepsilon)$.
  Then, by \cite[Lemma 11.1]{michael:56b}, there exists a continuous
  map $q:\sphere^{k}\to S$ such that $\rho(q(s),p(s))<\gamma(\xi)$ for
  every $s\in \sphere^{k}$. Hence, for every $s_1,s_2\in \sphere^k$,
  we have
  \begin{eqnarray*}
  d(q(s_1),q(s_2))& \le & d(q(s_1),p(s_1)) + d(p(s_1),p(s_2)) +
  d(p(s_2),q(s_2))\\
  & \le & \rho(q(s_1),p(s_1)) + d(p(s_1),p(s_2)) +
  \rho(p(s_2),q(s_2))\\
  & < & \gamma(\xi) + \eta(\varepsilon) + \gamma(\xi)\\
  & \le & \delta(\varepsilon/2)/4 + \delta(\varepsilon/2)/2
  +\delta(\varepsilon/2)/4=\delta(\varepsilon/2).
  \end{eqnarray*}
  So, $\diam_d(q(\sphere^k))<\delta(\varepsilon/2)$ and therefore, by
  hypothesis, there is a homotopy $h_{1}$ of $q$ to a constant map
  over a subset of $S$ of $d$-diameter$\ <\varepsilon/2$.  Let
  $\varphi(s,t)=Y$ for every $(s,t)\in\sphere^{k}\times \uin$, where
  $\uin=[0,1]$ denotes the closed unit interval.  Also, let
  $g(s,t)=q(s)$ for every $(s,t)\in \sphere^{k}\times \uin$, while
  $\ell:\sphere^{k}\times\{0,1\}$ be defined by
  $\ell|\sphere^{k}\times\{0\}=p$ and $\ell|\sphere^{k}\times\{1\}=q$.
  Since $\rho(\ell(x),g(x))<\gamma(\xi)$ for every
  $x\in\sphere^{k}\times\{0,1\}$, by \cite[Corollary
  4.2]{michael:56b}, $\ell$ can be extended to a continuous selection
  $h_{2}:\sphere^{k}\times \uin\to Y$ for $\varphi$ such that
  $\rho(h_{2}(x),g(x))<\xi\le\mu$ for every
  $x\in\sphere^{k}\times\uin$.  Thus, we get a homotopy $h_{2}$
  between $p$ and $q$, with $h_2(\sphere^k\times \uin)\subset
  B^\rho_\mu(S)$.  Finally, define $h$ to be the homotopy obtained by
  combining $h_{1}$ with $h_{2}$. Then, $h$ is a homotopy of $p$ with
  a constant map over a subset of $B_{\mu}^{\rho}(S)$ of
  $d$-diameter$\ <\varepsilon/2 + 2\xi<\varepsilon$.
\end{proof}

For a locally finite simplicial complex $M$ we use $|M|$ to denote the
\emph{polytope} on $M$, and $M^{k}$ to denote the \emph{$k$-skeleton}
of $M$. Also, for a locally finite cover ${\cal U}$ of $X$ we denote
by ${\cal N}({\cal U})$ the \emph{nerve} of ${\cal U}$, i.e. the
simplicial complex ${\cal N}({\cal U}) = \{\sigma \subset {\cal
  U}:\bigcap\sigma \neq \emptyset\}$.\medskip

Repeating precisely the proof of \cite[Lemma 6.1]{michael:56b} but now
using Lemma \ref{lm:functions} instead of \cite[Lemma
11.2]{michael:56b}, we get the following generalization of \cite[Lemma
6.1]{michael:56b}.

\begin{lemma}\label{lm:nerve-av}
  Let $E$ be a completely metrizable space, $Y\subset E$ be a
  $G_{\delta}$-set, $d\in \mathcal{D}(E)$, and let $\mathcal{S}\subset
  2^{Y}$ be $d$-uniformly equi-$LC^{n}$.  Also, let $\rho\in
  \mathcal{D}(Y)$ be a complete metric such that $\rho\ge d|Y\times Y$
  and $\mathcal{S}'=\mathcal{S}\cup\{Y\}$ is $\rho$-uniformly
  equi-$LC^{n}$. Then, to every $\varepsilon>0$ there corresponds
  $\beta(\varepsilon)>0$ with the following property\/$:$ If $X$ is a
  paracompact space, $\varphi:X\to \mathcal{S}$ is l.s.c., $g:X\to E$
  is a continuous map, with $g(x)\in
  B_{\beta(\varepsilon)}^{d}(\varphi(x))$ for every $x\in X$, and if
  $\mu>0$, then there exists a locally finite open cover $\mathcal{U}$
  of $X$ and a continuous $u:|\mathcal{N}^{n+1}(\mathcal{U})|\to Y$
  such that
  \[
  u(|\sigma|)\subset B_{\mu}^{\rho}(\varphi(x))\cap
  B_{\varepsilon}^{d}(g(x)),\ \hbox{for every $x\in\bigcap\sigma$ and
    $\sigma\in \mathcal{N}^{n+1}(\mathcal{U})$}.
  \]
\end{lemma}

We are now ready to prove the promised generalization of \cite[Theorem
4.1]{michael:56b}.

\begin{theorem}\label{th:avoid-open}
  Let $M$ be a completely metrizable space, $U\subset M$ be an open
  set, $d\in \mathcal{D}(M)$, and let $\mathcal{S}\subset
  \mathcal{F}(U)$ be $d$-uniformly equi-$LC^n$.  Then, for every
  $\varepsilon>0$ there exists $\nu(\varepsilon)>0$ with the following
  property\/$:$ If $X$ is a paracompact space, with $\dim(X)\le n+1$,
  $\varphi:X\to \mathcal{S}$ is l.s.c., and $g:X\to M$ is continuous,
  with $g(x)\in B_{\nu(\varepsilon)}^d(\varphi(x))$ for every $x\in
  X$, then $\varphi$ has a continuous selection $f$ such that
  $d(f(x),g(x))<\varepsilon$ for every $x\in X$.
\end{theorem}

\begin{proof}
  The proof involves arguments similar to those in \cite{michael:56b}.
  Namely, embed $(M,d)$ isometrically in a Banach space $(E,d)$, where
  $d$ is the metric on $E$ generated by the norm of $E$. Next, take an
  open set $Y\subset E$ such that $Y\cap M=U$.  Then, from one hand,
  $\mathcal{S}$ is $d$-uniformly equi-$LC^{n}$ as a family of subsets
  of $Y$. From another hand, $\mathcal{S}'=\mathcal{S}\cup\{Y\}$ is
  equi-$LC^{n}$ in $Y$ because $Y$ is an open subset of a Banach
  space. Since $Y$ is completely metrizable, there exists a complete
  metric $\tilde{d}\in \mathcal{D}(Y)$ such that $\tilde{d}\ge
  d|Y\times Y$. Therefore, by \cite[Proposition 2.1]{michael:56b},
  there also exists a metric $\rho\in \mathcal{D}(Y)$ such that
  $\rho\ge \tilde{d}$ and $\mathcal{S}'$ is $\rho$-uniformly
  equi-$LC^{n}$. In particular, $\rho$ is a complete metric in $Y$ and
  $\rho\ge d|Y\times Y$. Let $\gamma(\varepsilon)$ be as in
  \cite[Theorem 4.1]{michael:56b} applied to $\mathcal{S}'$ as a
  family which is $\rho$-uniformly equi-$LC^{n}$ in $Y$. Also, let
  $\beta(\varepsilon)$ be as in Lemma \ref{lm:nerve-av} applied to
  $d\in \mathcal{D}(E)$, $\mathcal{S}\subset 2^{Y}$, and $\rho\in
  \mathcal{D}(Y)$. Then, take $\nu(\varepsilon)=\beta(\varepsilon/2)$,
  and let us check that this works. So, let $\varphi:X\to \mathcal{S}$
  be l.s.c., and let $g:X\to M$ be continuous such that $g(x)\in
  B_{\beta(\varepsilon/2)}^d(\varphi(x))$ for every $x\in X$. Note
  that $g$ as a map from $X$ to $E$ has the same properties because
  $(M,d)$ is embedded isometrically in $(E,d)$. Hence, by Lemma
  \ref{lm:nerve-av}, applied with $\mu=\gamma(\varepsilon/2)$, there
  now exists a locally finite open cover $\mathcal{U}$ of $X$ and a
  continuous $u:|\mathcal{N}^{n+1}(\mathcal{U})|\to Y$ such that
  \begin{equation}
    \label{eq:ord-1}
  u(|\sigma|)\subset B_{\gamma(\varepsilon/2)}^{\rho}(\varphi(x))\cap
  B_{\varepsilon/2}^{d}(g(x)),\ \hbox{for every $x\in\bigcap\sigma$ and
  $\sigma\in \mathcal{N}^{n+1}(\mathcal{U})$}.
  \end{equation}
  Since $\dim(X)\le n+1$, there also exists an open cover
  $\{V_{U}:U\in\mathcal{U}\}$ of $X$ such that $V_{U}\subset U$ for
  every $U\in \mathcal{U}$, and
  \begin{equation}
    \label{eq:ord-2}
    \{U\in \mathcal{U}: x\in
  V_{U}\}\in\mathcal{N}^{n+1}(\mathcal{U}),\quad \hbox{for every $x\in X$.}
  \end{equation}
  Finally, take a partition of unity $\{\xi_U:U\in \mathcal{U}\}$
  index-subordinated to the cover $\{V_{U}:U\in\mathcal{U}\}$ of $X$.
  Then, by (\ref{eq:ord-2}), $k=u\circ\xi:X\to Y$ defines a continuous
  map, where $\xi:X\to \left|\mathcal{N}(\mathcal{U})\right|$ is the
  canonical map $\xi(x)=\sum\{\xi_{U}(x)\cdot U:U\in\mathcal{U}\}$,
  $x\in X$. According to (\ref{eq:ord-1}) and (\ref{eq:ord-2}), this
  map has also the property that
  \[
  k(x)\in B_{\gamma(\varepsilon/2)}^{\rho}(\varphi(x))\cap
  B_{\varepsilon/2}^{d}(g(x)),\quad \hbox{for every $x\in X$.}
  \]
  Then, by \cite[Theorem 4.1]{michael:56b}, there exists a continuous
  selection $f$ for $\varphi$ such that
  $\rho(f(x),k(x))<\varepsilon/2$ for every $x\in X$. Since $\rho\ge
  d|Y\times Y$, this finally implies that, for every $x\in X$,
  \begin{eqnarray*}
  d(f(x),g(x))&\le&
  d(f(x),k(x))+ d(k(x),g(x))\\
  &\le& \rho(f(x),k(x))+
  d(k(x),g(x))\\
  &<&\varepsilon/2 +\varepsilon/2=\varepsilon,
  \end{eqnarray*}
  which completes the proof.
\end{proof}

Theorem \ref{th:avoid-open} has some interesting consequences. One of
them is related to the requirement in \cite[Theorem 4.1]{michael:56b}
the family $\mathcal{S}\subset \mathcal{F}(Y)$ to be $d$-uniformly
equi-$LC^{n}$ with respect to a complete metric $d\in\mathcal{D}(Y)$.
Namely, by taking $U=M$ in Theorem \ref{th:avoid-open}, we have the
following immediate result.

\begin{corollary}\label{cr:completeness}
  Let $X$ be a paracompact space, with $\dim(X)\le n+1$, $Y$ be a
  completely metrizable space, $d\in \mathcal{D}(Y)$, and let
  $\mathcal{S}\subset \mathcal{F}(Y)$ be $d$-uniformly equi-$LC^n$.
  Then, for every $\varepsilon>0$ there exists $\nu(\varepsilon)>0$
  with the following property: If $\varphi:X\to \mathcal{S}$ is
  l.s.c.\ and $g:X\to Y$ is continuous, with $g(x)\in
  B_{\nu(\varepsilon)}^d(\varphi(x))$ for every $x\in X$, then
  $\varphi$ has a continuous selection $f$ such that
  $d(f(x),g(x))<\varepsilon$ for every $x\in X$.
\end{corollary}

In our next consequence, a function $\mu:X\to \R$ is
\emph{lower semi-continuous} if the set $\{x\in X:\mu(x)>r\}$
is open in $X$ for every $r\in \R$.

\begin{corollary}\label{th:fine-app}
  Let $M$ be a completely metrizable space, $U\subset M$ be an open
  set, $d\in \mathcal{D}(M)$, and let $\mathcal{T}\subset
  \mathcal{F}(U)$ be $d$-uniformly equi-$LC^n$.  Then, for every
  paracompact space $X$, with $\dim(X)\le n+1$, and every lower
  semi-continuous function $\mu:X\to (0,+\infty)$ there exists a
  continuous function $\alpha:X\to (0,+\infty)$ with the following
  property\/$:$ If $\varphi:X\to \mathcal{T}$ is l.s.c.\ and $g:X\to
  M$ is continuous, with \hbox{$g(x)\in
    B_{\alpha(x)}^{d}(\varphi(x))$} for every $x\in X$, then there
  exists a continuous selection $f$ for $\varphi$ such that
  $d(f(x),g(x))<\mu(x)$ for every $x\in X$.
\end{corollary}

\begin{proof}
  The proof follows an idea in \cite[Proof that Theorem 4.1 implies
  Theorem 1.3]{michael:56b}. Namely, let $X$ and $\mu$ be as in the
  hypotheses. Then,
  \begin{equation}\label{eq:lsc-inc}
  V_i=\{x\in X: \mu(x)> 1/i\},\quad i\in \N,
  \end{equation}
  defines an open increasing cover of $X$ because $\mu$ is lower
  semi-continuous. Since $X$ is paracompact (hence, normal and
  countably paracompact as well), there exists another open cover
  $\{W_i:i\in \N\}$ of $X$ such that
   \begin{equation}
  \label{eq:inc-cov}
  \overline{W_i}\subset
  V_i\cap W_{i+1},\quad \hbox{for every $i\in \N$.}
  \end{equation}
  Next, for every $\varepsilon>0$, let $\nu(\varepsilon)\le
  \varepsilon$ be as in Theorem \ref{th:avoid-open} applied to the
  family $\mathcal{S}=\mathcal{T}\cup\{\{y\}:y\in U\}$.  Define a
  decreasing function $\eta:\N\to (0,+\infty)$ by letting for $i\in
  \N$ that
  \begin{equation}
    \label{eq:decreas-fun}
    \eta(i)=\min\left\{\nu\left(\nu\left(\frac 1j\right )\cdot
     \nu\left(\frac{1}{j+1}\right)\right): 1\le j\le i\right\}.
  \end{equation}
  Now, for every $x\in X$, let $i(x)=\min\{i\in\N: x\in W_i\}$. Then,
  define a function $\beta:X\to (0,+\infty)$ by
  \begin{equation}
    \label{eq:n-rest}
    \beta(x)=\eta(i(x)),\quad x\in X.
  \end{equation}
  Note that $\beta$ is lower semi-continuous. Indeed, if $x\in
  W_{i(x_{0})}$ for some point $x_{0}\in X$, then $i(x)\le i(x_{0})$
  implies $\beta(x)=\eta(i(x))\ge \eta(i(x_{0}))=\beta(x_{0})$.  Since
  $X$ is paracompact, by a result of \cite{dieudonne:44} (see, also,
  \cite{dowker:51,katetov:51}), there exists a continuous function
  $\alpha:X\to (0,+\infty)$ such that $\alpha(x)\le\beta(x)$ for every
  $x\in X$. This $\alpha$ is as required. Indeed, take an l.s.c.\
  $\varphi:X\to \mathcal{T}$ and a continuous $g:X\to M$ with
  $d(g(x),\varphi(x))<\alpha(x)\le\beta(x)$ for every $x\in X$. Also,
  let $A_0=\emptyset$, and $A_i=\overline{W_i}\setminus W_i$ for every
  $i\in \N$. Then, each $A_i$ is a paracompact space as a closed
  subset of $X$, and $\dim(A_i)\le n+1$. On the other hand, $x\in
  A_{i}$ implies $x\notin W_{i}$, so $i(x)>i$. Therefore, by
  (\ref{eq:decreas-fun}) and (\ref{eq:n-rest}), we now get that
  \[
  d(g(x),\varphi(x))<\beta(x)=\eta(i(x))\le \eta(i)\le
  \nu(\nu(1/i)\cdot \nu(1/(i+1))).
  \]
  Hence, by Theorem \ref{th:avoid-open}, each $\varphi|A_i$ has a
  continuous selection $h_i:A_i\to U$ with
  \begin{equation}
  \label{eq:fin-app}
  d(h_i(x),g(x))<\nu(1/i)\cdot \nu(1/(i+1)),\quad
  \hbox{for every $x\in A_i$ and $i\in\N$}.
  \end{equation}
  Let $B_i=A_{i-1}\cup A_i$ and $X_i=\overline{W_{i}}\setminus
  W_{i-1}$ for every $i\in\N$, where $W_0=\emptyset$. Next, define
  $k_i:B_i\to U$ by $k|A_j=h_j$, $j=i-1,i$. Finally, define
  $\varphi:X_{i}\to \mathcal{S}$ by $\varphi(x)=\{k_{i}(x)\}$ if $x\in
  B_{i}$ and $\varphi_{i}(x)=\varphi(x)$ otherwise. By \cite[Example
  1.3$^{*}$]{michael:56a}, each $\varphi_{i}$ is l.s.c. Then, by
  \eqref{eq:fin-app} and Theorem \ref{th:avoid-open}, each $\varphi_i$
  has a continuous selection $f_i:X_i\to U$ such that
  $d(f_i(x),g(x))<1/i$ for every $x\in X_i$. In fact, $f_{i}$ is a
  continuous extension of $k_{i}$. Since, by (\ref{eq:inc-cov}),
  $X_i\cap X_{i+1}=A_{i}$, while $f_i|A_i=h_i=f_{i+1}|A_{i}$, we may
  now define $f:X\to Y$ by $f|X_i=f_i$. Clearly, $f$ is a continuous
  selection for $\varphi$ and, by \eqref{eq:lsc-inc} and
  (\ref{eq:inc-cov}), $x\in X_i\subset V_i$ implies
  $d(f(x),g(x))=d(f_i(x),g(x))<1/i<\mu(x)$.
\end{proof}

\section{A construction of set-valued mappings}

In this section, we continue the preparation for the proof of Theorem
\ref{th:fdver} involving a construction of \cite{michael:89}. Namely,
let $(Y,d)$ be a metric space, $\psi:X\to \mathcal{P}(Y)$ be a mapping
with a closed graph, and let $M_{\psi}=(X\times Y)\setminus
\graph(\psi)$. Then, for every $(x,y)\in M_{\psi}$, we consider the
set $\Delta(x,y)$ of all $\delta\in (0,1]$ for which there exists a
neighbourhood $U_{\delta}$ of $x$ such that $U_{\delta}\times
B^{d}_{\delta}(y)\subset M_{\psi}$. Note that
$\Delta(x,y)\ne\emptyset$ for every $(x,y)\in M_{\psi}$. Now, we
define a map $u:X\times Y\to \uin$ by letting for $(x,y)\in X\times Y$
that
\begin{equation}
  \label{eq:def-u}
 u(x,y)=\begin{cases}
    \sup\Delta(x,y), & \hbox{if}\ \ (x,y)\in M_{\psi},\\
    0, & \hbox{otherwise}.
              \end{cases}  
\end{equation}

The following properties of $u$ were actually demonstrated in
\cite[Lemma 3.1]{michael:89}. 

\begin{proposition}[\cite{michael:89}]\label{pr:michael-mod1}
  Let $(Y,d)$ be a metric space, $\psi:X\to \mathcal{P}(Y)$ be a
  mapping with a closed graph, and let $u:X\times Y\to \uin$ be
  defined as in $(\ref{eq:def-u})$. Then, 
  \begin{enumerate}
  \item[\emph{(a)}] $u$ is a lower semi-continuous function.
  \item[\emph{(b)}] $|u(x,y_{1})-u(x,y_{2})|\le d(y_{1},y_{2})$ for
    every $x\in X$ and $y_{1},y_{2}\in Y$.  
  \end{enumerate}
\end{proposition}

\begin{proof}
  Our situation is slightly different from that in \cite[Lemma
  3.1]{michael:89}, so we provide some arguments to the proof. Since
  $u(x,y)\ge 0$ for every $(x,y)\in X\times Y$, to show (a), it
  suffices to consider only the points of $M_{\psi}$. To this end, we
  repeat precisely the arguments in \cite[Lemma
  2.1]{gutev-valov:00}. Namely, take a point $(x_{0},y_{0})\in
  M_{\psi}$ and $a\in \uin$, with $a<u(x_{0},y_{0})$. Then, there
  exists $\delta_{0}\in (a,u(x_{0},y_{0})]$ and a neighbourhood
  $U_{0}$ of $x_{0}$ such that $U_{0}\times
  B_{\delta_{0}}^{d}(y_{0})\subset M_{\psi}$. Let
  $\rho=(\delta_{0}-a)/2$ and $\delta=a+\rho$. It now follows
  that $\delta\in \Delta(x,y)$ for every $(x,y)\in U_{0}\times
  B_{\rho}^{d}(y_{0})$ because 
  \[
  U_{0}\times B_{\delta}^{d}(y)\subset U_{0}\times
  B_{\delta+\rho}^{d}(y_{0}) \subset U_{0}\times
  B_{\delta_{0}}^{d}(y_{0})\subset M_{\psi}, \quad (x,y)\in
  U_{0}\times B_{\rho}^{d}(y_{0}).
  \]
  Hence, in particular, $u(x,y)\ge \delta>a$ for every $(x,y)\in U_{0}\times
  B_{\rho}^{d}(y_{0})$.\smallskip

  The statement in (b) is the same as in \cite[Lemma
  3.1]{michael:89}. 
\end{proof}

We are now ready to state the main result of this section which is an
analogue to \cite[Lemma 3.3]{michael:89}. 

\begin{theorem}\label{th:const-1}
  Let $X$ be a space, $(Y,d)$ be a complete metric space, $\varphi
  :X\rightarrow \mathcal{F}(Y)$ be l.s.c., with $\{\varphi (x): x\in
  X\}$ $d$-uniformly equi-$LC^{n-1}$, and let $\psi: X\rightarrow
  \mathcal{P}(Y)$ be a mapping with a closed graph such that
  $\psi(x)\cap\varphi(x)$ is a $\sigma Z_{n}$-set in $\varphi(x)$ for
  every $x\in X$. Define another set-valued mapping $\Phi:X\to
  2^{Y\times (0,1]}$ by
  \[
  \Phi(x)=\{(y,t)\in Y\times (0,1]: y\in \varphi(x)\setminus \psi(x)\
  \hbox{and}\ t\in (0,u(x,y)]\},\quad x\in X.
  \]
  Then$:$
  \begin{enumerate}
  \item[\emph{(a)}] $\Phi$ is l.s.c.
  \item[\emph{(b)}] $\Phi(x)$ is a non-empty closed subset of $Y\times
    (0,1]$, for every $x\in X$.
  \item[\emph{(c)}] The family $\{\Phi(x):x\in X\}$ is $d\times
    e$-uniformly equi-$LC^{n-1}$, where $d\times e$ is the box metric
    $d\times e((y_{1},t_{1}),(y_{2},t_{2}))=\max\{d(y_{1},y_{2}),
    e(t_{1},t_{2})\}$ on $Y\times (0,1]$ generated by $d$ and the
    Euclidean metric $e$ on $(0,1]$.
  \end{enumerate}
\end{theorem}

To prepare for the proof of Theorem \ref{th:const-1} we need the
following simple observations about $Z_{n}$-sets in complete metric
spaces, the first of which is well-known.

\begin{lemma}\label{pr:sigzn}
  Let $Y$ be a completely metrizable space, and let $T\subset Y$ be
  closed. Then, $T$ is a $\sigma Z_{n}$-set in $Y$ if and only if it
  is a $Z_{n}$-set in $Y$.
\end{lemma}

\begin{proposition}\label{pr-unif-z}
  Let $(Y,d)$ be a complete metric space, $\mathcal{S}\subset
  \mathcal{F}(Y)$ be $d$-uniformly equi-$LC^{n-1}$, and let
  $\{Z_S\subset S:S\in \mathcal{S}\}$ be a family of closed sets in
  $Y$ such that each $Z_S$ is a $ Z_{n}$-set in $S$, $S\in
  \mathcal{S}$.  Then, $\{S\setminus Z_S: S\in \mathcal{S}\}$ is also
  $d$-uniformly equi-$LC^{n-1}$.
\end{proposition}

\begin{proof}
  In case $n=0$ the proof is trivial. So, we suppose that $n\ge 1$.
  Let $\delta(\varepsilon)\le \varepsilon$ be as in the definition of
  $d$-uniformly equi-$LC^{n-1}$ of $\mathcal{S}$. We are going to show
  that $\{S\setminus Z_S: S\in \mathcal{S}\}$ is $d$-uniformly
  equi-$LC^{n-1}$ with respect to $\delta(\varepsilon/3)$. Take
  $\varepsilon>0$, $S\in \mathcal{S}$, and a continuous map
  $p:\sphere^{m}\to S\setminus Z_S$ (for some $0\le m\le n-1$) such
  that $\diam_{d}(p(\sphere^{m}))\le \delta(\varepsilon/3)$. Let
  \[
  \cone (\sphere^{m})=\{(1-t)\cdot \vartheta+t\cdot s: (s,t)\in
  \sphere^{m}\times \uin\}
  \]
  be the \emph{cone} of $\sphere^{m}$ with a vertex $\vartheta$. Also,
  take an $a\in (0,1)$, and consider the copy
  \[
  \cone_{a} (\sphere^{m})=\{(1-t)\cdot\vartheta+t\cdot s: (s,t)\in
  \sphere^{m}\times [0,a]\}
  \]
  of the cone $\cone(\sphere^{m})$. Whenever $t>0$, for convenience,
  let $\sphere^{m}_{t}=\sphere^{m}\times\{t\}$ be the corresponding
  copy of $\sphere^{m}$ in $\cone(\sphere^{m})$. Then, according to
  the properties of $\delta(\varepsilon/3)$, there exists a continuous
  $q:\cone_{a}(\sphere^m)\to S$ such that
  \begin{equation}
    \label{eq:ext-1}
    q|\sphere^m_{a}=p
    \quad\hbox{and}\quad
    \diam_{d}(q(\cone_{a}(\sphere^m)))<\frac{\varepsilon}{3}.
  \end{equation}
  Let $\varepsilon_{1}\in (0,\varepsilon/6]$ be such that
  $B_{2\varepsilon_{1}}^{d}(p(\sphere^{m}))\cap Z_{S}=\emptyset$.
  Next, let $\gamma(\varepsilon_1)\le \varepsilon_1$ be as in
  \cite[Corollary 4.2]{michael:56b} applied to the family
  $\mathcal{S}$.  Since $Z_S$ is a $Z_n$-set in $S$, the set
  $C(\cone_{a}(\sphere^m), S\setminus Z_S)$ is dense in
  $C(\cone_{a}(\sphere^m), S)$ with respect to the uniform topology.
  Therefore, there exists a continuous $\ell:\cone_{a}(\sphere^m)\to
  S\setminus Z_S$ such that
    \begin{equation}
      \label{eq:apr-1}
      d(\ell(x),q(x))<\gamma(\varepsilon_1),\ x\in \cone_{a}(\sphere^m).
    \end{equation}
    For later use, let us observe that, by (\ref{eq:ext-1}), this
    implies
  \begin{equation}
    \label{eq:diam-2}
    \diam_{d}(\ell(\cone_{a}(\sphere^{m})))<2\gamma(\varepsilon_{1})
    +\diam_{d}(q(\cone_{a}(\sphere^m)))\le \frac{2\varepsilon}{6}
    +\frac{\varepsilon}{3}=\frac{2\varepsilon}{3}.
  \end{equation}
  Define now a map $k:\cone(\sphere^m)\to S$ by letting for every
  $(s,t)\in \sphere^{m}\times\uin$ that
  \[
  k((1-t)\cdot \vartheta +t\cdot s)=\begin{cases}
    \ell((1-t)\cdot\vartheta +t\cdot s), & \hbox{if}\ \ t\le a,\\
    \ell((1-a)\cdot\vartheta +a\cdot s), & \hbox{otherwise}.
              \end{cases}
  \]
  Also, define another map $g:\cone_{a}(\sphere^{m})\cup
  \sphere^{m}_{1}\to S$ by $g|\cone_{a}(\sphere^{m})=\ell$, and
  $g|\sphere^{m}_{1}=p$. Then, by (\ref{eq:ext-1}) and
  (\ref{eq:apr-1}), $d(g(x),k(x))<\gamma(\varepsilon_{1})$ for every
  $x\in \cone_{a}(\sphere^{m})\cup \sphere^{m}_{1}$. So, by
  \cite[Corollary 4.2]{michael:56b} (with $\varphi(x)=S$ for every
  $x\in \cone(\sphere^{m})$), there exists a continuous
  $h:\cone(\sphere^m)\to S$ such that
    \begin{equation}
      \label{eq:ext-f}
      h|\cone_{a}(\sphere^m)=\ell,\quad h|\sphere^m_{1}=p, \quad
      \hbox{and}\quad
      d(h(x),k(x))<\varepsilon_1,\ x\in \cone(\sphere^m).
    \end{equation}
    Clearly, $h$ is a homotopy connecting $p$ and a constant map.
    Moreover, by (\ref{eq:diam-2}) and (\ref{eq:ext-f}),
    \[
    \diam_{d}(h(\cone(\sphere^{m})))\le 2\varepsilon_{1}
    +\diam_{d}(\ell(\cone_{a}(\sphere^{m})))< \frac{2\varepsilon}{6} +
    \frac{2\varepsilon}{3} = \varepsilon.
    \]
    On the other hand, by (\ref{eq:ext-1}) and (\ref{eq:apr-1}),
    \begin{eqnarray*}
    h(\cone(\sphere^{m})\setminus \cone_{a}(\sphere^{m}))&\subset&
    B_{\varepsilon_{1}}^{d}(k(\cone(\sphere^{m})\setminus \cone_{a}(\sphere^{m})))\\
    &=& B_{\varepsilon_{1}}^{d}(\ell(\sphere^{m}_{a}))\\
     &\subset&
    B_{\varepsilon_{1}+\gamma(\varepsilon_{1})}^{d}(p(\sphere^{m}))\\
     &\subset&
    B_{2\varepsilon_{1}}^{d}(p(\sphere^{m})) \subset S\setminus Z_{S}.
    \end{eqnarray*}
    Since, by (\ref{eq:ext-f}), $h(\cone_{a}(\sphere^{m}))=
    \ell(\cone_{a}(\sphere^{m})) \subset S\setminus Z_{S}$, this
    finally implies that $h(\cone(\sphere^{m}))\subset S\setminus
    Z_{S}$ which completes the proof.
\end{proof}

\begin{proof*}{Proof of Theorem \ref{th:const-1}}
  Let $X$, $(Y,d)$, $\varphi$ and $\psi$ be as in that theorem.  The
  statements of (a) and (b) were established in \cite[Lemma
  3.3]{michael:89}, so we prove only (c). According to Lemma
  \ref{pr:sigzn} and Proposition \ref{pr-unif-z}, the family
  $\mathcal{T}=\{\varphi(x)\setminus \psi(x): x\in X\}$ is
  $d$-uniformly equi-$LC^{n-1}$. Then, let $\delta(\varepsilon)$ be as
  in the definition of $d$-uniformly equi-$LC^{n-1}$ of $\mathcal{T}$.
  We are going to show that $\{\Phi(x):x\in X\}$ is $d\times
  e$-uniformly equi-$LC^{n-1}$ with respect to
  $\delta(\varepsilon/2)$. So, take an $\varepsilon>0$, a point $x\in
  X$, and a continuous map $p:\sphere^{k}\to \Phi(x)$ for some $k$,
  with $0\le k\le n-1$, such that $\diam_{d\times
    e}(p(\sphere^{k}))<\delta(\varepsilon/2)$. Also, let
  $\pi_{1}:Y\times (0,1]\to Y$ and $\pi_{2}:Y\times (0,1]\to (0,1]$ be
  the corresponding projections. Then, $p_{1}=\pi_{1}\circ
  p:\sphere^{k}\to \varphi(x)\setminus \psi(x)$ is a continuous map,
  with $\diam_{d}(p_{1}(\sphere^{k}))<\delta(\varepsilon/2)$. Hence,
  $p_{1}$ can be extended to a continuous map $q_{1}:\ball^{k+1}\to
  \varphi(x)\setminus \psi(x)$ such that
   $\mu=\diam_{d}(q_{1}(\ball^{k+1}))<\varepsilon/2<\varepsilon$. In the
  same way, $p_{2}=\pi_{2}\circ p:\sphere^{k}\to (0,1]$ is a
  continuous map, with
  $\diam_{e}(p_{2}(\sphere^{n-1}))<\delta(\varepsilon/2)$. According
  to the definition of $\Phi(x)$, we now have that $p_{2}(s)\le
  u(x,p_{1}(s))$ for every $s\in \sphere^{k}$. Fix an arbitrary point
  $s_{0}\in \sphere^{k}$. Then, by Proposition \ref{pr:michael-mod1},
  we have that
  \[
  |u(x,q_{1}(b))-u(x,q_{1}(s_{0}))|\le
   d(q_{1}(b),q_{1}(s_{0}))\le \mu <\frac{\varepsilon}{2}, \quad
   b\in\ball^{k+1}.
  \]
  Hence, $B_{\mu}^{e}(p_{2}(s_{0}))\cap (0,u(x,q_{1}(b))]\ne
  \emptyset$, $b\in \ball^{k+1}$, because $p_{2}(s_{0})\le
  u(x,p_{1}(s_{0}))=u(x,q_{1}(s_{0}))$. Since, by Proposition
  \ref{pr:michael-mod1}, $u(x,q_{1}(b))$, $b\in \ball^{k+1}$, is a
  continuous map, we may now define an l.s.c.\ mapping
  $\theta:\ball^{k+1}\to \mathcal{F}((0,1])$ by
  \[
  \theta(b)=\overline{B_{\mu}^{e}(p_{2}(s_{0}))\cap
  (0,u(x,q_{1}(b))]}^{\,(0,1]}, \quad b\in \ball^{k+1}.
  \]
  Then, by \cite[Theorem
  3.2$''$]{michael:56a}, $\theta$ has a continuous selection
  $q_{2}:\ball^{k+1}\to (0,1]$ because $\theta$ is convex-valued. Note
  that $q_{2}$ is an extension of $p_{2}$ such that
  $e(q_{2}(b),q_{2}(s_{0}))\le \mu$ for every
  $b\in\ball^{k+1}$, i.e.\ $\diam_{e}(q_{2}(\ball^{k+1}))\le 2\mu
  <\varepsilon$. Then, the diagonal map $q=q_{1}\Delta
  q_{2}:\ball^{k+1}\to \Phi(x)$ is an extension of $p$, with
  $\diam_{d\times e}(q(\ball^{k+1}))<\varepsilon$. 
\end{proof*}

\section{Proof of Theorem \ref{th:fdver}}

We finalize the preparation for the proof of Theorem \ref{th:fdver}
with the following result which may have an independent
interest.

\begin{theorem}\label{th:st-1}
  Let $X$ be paracompact, with $\dim (X)\leq n$, $Y$ be a completely
  metrizable space, $\varphi :X\rightarrow \mathcal{F}(Y)$ be l.s.c.\
  such that $\{\varphi (x): x\in X\}$ is equi-$LC^{n-1}$, and let $\psi:
  X\rightarrow \mathcal{P}(Y)$ be a mapping with a closed graph
  such that $\psi(x)\cap\varphi(x)$ is a $\sigma Z_{n}$-set in
  $\varphi(x)$ for every $x\in X$. Also, let $f$ be a continuous
  selection for $\varphi$, and let $\mu\in C(X,(0,+\infty))$. Then for
  every $\rho\in \mathcal{D}(Y)$, the mapping $\varphi$ has a
  continuous selection $g$ such that $g(x)\notin \psi(x)$ and
  $\rho(g(x),f(x))<\mu(x)$ for every $x\in X$.
\end{theorem}

\begin{proof}
  Let $\rho\in \mathcal{D}(Y)$, and let $\tilde{d}\in \mathcal{D}(Y)$
  be a complete metric, with $\tilde{d}\ge \rho$. Then, by
  \cite[Proposition 2.1]{michael:56b} (see, also, \cite[Theorem
  1]{dugundji-michael:56}), there exists another complete metric $d\in
  \mathcal{D}(Y)$ such that $\mathcal{S}=\{\varphi(x):x\in X\}$ is
  $d$-uniformly equi-$LC^{n-1}$ and $d\ge \tilde{d}\ge \rho$.  Now,
  let $u:X\times Y\to \uin$ be defined as in (\ref{eq:def-u}) using
  the metric $d$, and let $\Phi:X\to \mathcal{F}(Y\times (0,1])$ be
  defined as in Theorem \ref{th:const-1}. By Proposition
  \ref{pr:michael-mod1}, the function $r:X\to \uin$, defined by
  $r(x)=u(x,f(x))$, $x\in X$, is lower semi-continuous. Since $X$ is
  normal, there now exists a continuous function $v:X\to \uin$ such
  that $v(x)\le r(x)$, $x\in X$, see \cite{engelking:89}. Note that
  the diagonal map $\ell=f\Delta v:X\to Y\times\uin$ is a selection
  for the closure of $\Phi$, i.e.\ $\ell(x)\in\overline{\Phi(x)}$ for
  every $x\in X$. Thus, in particular, $\ell(x)\in
  B_{\alpha(x)}^{d\times e}(\Phi(x))$, $x\in X$, for every $\alpha\in
  C(X,(0,+\infty))$. On the other hand, by Theorem \ref{th:const-1},
  the family $\{\Phi(x):x\in X\}$ is $d\times e$-uniformly
  equi-$LC^{n-1}$.  Hence, by Corollary \ref{th:fine-app}, $\Phi$ has
  a continuous selection $p:X\to Y\times (0,1]$ such that $d\times
  e$-distance in $Y\times \uin$ between the points $p(x)$ and
  $\ell(x)$ is less than $\mu(x)$ for every $x\in X$ (recall that
  $Y\times (0,1]$ is open in $Y\times \uin$).  Then, $g=\pi_{1}\circ
  p:X\to Y$ is as required, where $\pi_{1}:Y\times (0,1]\to Y$ is the
  projection.
\end{proof}

By Theorem \ref{th:st-1} we have the following consequence.

\begin{corollary}\label{cr:fine-dense}
  Let $X$ be paracompact, with $\dim (X)\leq n$, $Y$ be a completely
  metrizable space, $\varphi :X\rightarrow \mathcal{F}(Y)$ be l.s.c.\ 
  such that $\{\varphi (x): x\in X\}$ is equi-$LC^{n-1}$ and each
  $\varphi(x)$ is $C^{n-1}$, $x\in X$, and let $\psi: X\rightarrow
  \mathcal{P}(Y)$ be a mapping with a closed graph such that
  $\psi(x)\cap\varphi(x)$ is a $\sigma Z_{n}$-set in $\varphi(x)$ for
  every $x\in X$. Then, the set $\sel(\varphi)$ of all continuous
  selections for $\varphi$, endowed with the fine topology, is a Baire
  space and $\mathcal{M}(\psi)\cap \sel(\varphi)=\{f\in \sel(\varphi):
  f(x)\notin \psi(x)\ \hbox{for every $x\in X$}\}$ is open and dense
  in $\sel(\varphi)$.
\end{corollary}

\begin{proof}
  Take a complete bounded metric $d\in \mathcal{D}(Y)$. Every subset
  of $C(X,Y)$ which is closed with respect to the uniform topology
  generated by $d$ is a Baire space in the fine topology, see
  \cite[Lemma 3.2]{michael:88a}.  Obviously, $\sel(\varphi)$ is
  uniformly closed in $C(X,Y)$ with respect to $d$, so it has the
  Baire property.  Let us show that
  $\mathcal{M}(\psi)\cap\sel(\varphi)$ is open in $\sel(\varphi)$
  (note that, by \cite[Theorem 1.2]{michael:56b} and Theorem
  \ref{th:st-1}, $\mathcal{M}(\psi)\cap \sel(\varphi)\ne\emptyset$).
  Take an $f\in\mathcal{M}(\psi)\cap \sel(\varphi)$, and let
  $u:X\times Y\to \uin$ be as in (\ref{eq:def-u}). Then, by
  Proposition \ref{pr:michael-mod1}, $\eta(x)=u(x,f(x))$, $x\in X$,
  defines a lower semi-continuous function $\eta:X\to (0,+\infty)$.
  Hence, by \cite{dieudonne:44,dowker:51,katetov:51}, there exists an
  $\alpha\in C(X,(0,+\infty))$, with $\alpha (x)<\eta(x)$ for every
  $x\in X$.  Therefore
  $V(f,\alpha)\cap\sel(\varphi)\subset\mathcal{M}(\psi)\cap\sel(\varphi)$,
  so $\mathcal{M}(\psi)\cap \sel(\varphi)$ is open in $\sel(\varphi)$.
  That $\mathcal{M}(\psi)\cap\sel(\varphi)\subset\sel(\varphi)$ is
  dense, it follows by Theorem \ref{th:st-1}.
\end{proof}

Now, we complete the proof of Theorem \ref{th:fdver} in the following
way. Let $X$, $Y$, $\varphi$ and $\psi$ be as in that theorem. Since
$\graph(\psi)$ is an $F_\sigma$-subset of $X\times Y$, there are
mappings $\psi_k:X\to \mathcal{P}(Y)$, $k\in\N$, such that each
$\psi_k$, $k\in\N$, has a closed graph and
$\psi(x)=\bigcup\{\psi_k(x):k\in\N\}$. Then, by Corollary
\ref{cr:fine-dense}, each set $\mathcal{M}(\psi_k)\cap \sel(\varphi)$,
$k\in\N$, is open and dense in $\sel(\varphi)$ with respect to the
fine topology, while $\sel(\varphi)$ is a Baire space.  Therefore,
$\mathcal{M}(\psi)\cap \sel(\varphi)=\bigcap\{\mathcal{M}(\psi_k)\cap
\sel(\varphi):k\in\N\}$ is a dense $G_\delta$-subset of
$\sel(\varphi)$ which completes the proof.\bigskip
 
\noindent\textbf{Remark.} It should be mentioned that Theorem
\ref{th:fdver} might be compare with \cite[Lemma 3.2]{michael:88a}
(see, also, \cite[Theorem 3.1]{michael:88a}) which works in the
special case of a convex-valued $\varphi$ and a continuous $\psi$. In
fact, one of the main goals of Theorem \ref{th:fdver} is to avoid the
assumption on $\varphi$ to be convex-valued. Namely, in case $Y$ is a
Banach space and $\varphi$ is convex-valued, Theorem \ref{th:fdver}
has a relative shorter proof based either on the technique developed
in \cite{michael:89} or on that stated in \cite{gutev-valov:00}.
Despite of that, in both cases the proof will rely on some results of
\cite{michael:56b}, hence will not make the corresponding arguments
simpler.  \bigskip

We complete this section showing that the following result of E.
Michael \cite[Theorem 5.2]{michael:89} can be derived from our Theorem
\ref{th:fdver}. In what follows, for a linear topological space
$Y$, we let $\mathcal{F}_c(Y)=\{S\in\mathcal{F}(Y): S\ \hbox{is
  convex}\}$.

\begin{theorem}\label{la}
  Let $X$ be a paracompact space, $Y$ be a Banach space, $\varphi
  :X\rightarrow\mathcal{F}_c(Y)$ be l.s.c., and let $\psi_n:X\to
  \mathcal{P}(Y)$, $n\in\mathbb{N}$, be a sequence of mappings such
  that each $\psi_n$ has a closed graph and $\dim (X)<\dim(\varphi
  (x))-\dim(\conv(\psi_n(x)\cap \varphi(x)))$ for every $x\in X$. Then
  $\varphi$ has a continuous selection $f$ such that
  $f(x)\not\in\bigcup\{\psi_n(x):n\in\N\}$ for every $x\in X$.
\end{theorem}

\begin{proof}
  Suppose that $\dim (X)=m$. To apply Theorem \ref{th:fdver}, it
  suffices to check that, in this situation, $\psi_n(x)\cap \varphi
  (x)$ is a $Z_m$-set in $\varphi (x)$ for all $x\in X$ and $n\in\N$.
  To this end, observe that, by \cite[Lemma 2.1]{michael:88a}, the
  sets $\varphi (x)\backslash\psi _n(x)$ are $C^{m-1}$ and, obviously,
  they are $LC^{m-1}$ because each $\varphi (x)$ is convex.  Moreover,
  if $V$ is an open ball in $\varphi (x)$, then
  \[
  m+\dim(\conv(\psi _n(x)\cap V))<\dim (V)\ \ \hbox{for every
    $n\in\N$.}
  \]
  Hence, again by \cite[Lemma 2.1]{michael:88a}, $V\backslash \psi
  _n(x)$ is $C^{m-1}$. Then, according to \cite[Theorem 2.8 and
  Corollary 3.3]{torunczyk:78}, each $\psi_n(x)\cap \varphi (x)$ is a
  $Z_m$-set in $\varphi (x)$.
\end{proof}

\section{Strongly countable-dimensional spaces and selections} 

In this section we provide some applications of Theorem \ref{th:fdver}
related to \emph{strongly countable-dimensional spaces} (i.e, spaces
which are a countable union of closed finite-dimensional subspaces).
Our first result is a possible selection analogue of a result
established by E.\ Pol \cite[Corollary 4.4]{epol:88b}.

\begin{theorem}\label{th:count-case}
  Let $X_n$, $n<\omega$, be closed subsets of a paracompact space $X$,
  with $\dim (X_n)\leq n$, $Y$ be a countable product of Banach
  spaces, $\varphi :X\rightarrow \mathcal{F}_c(Y)$ be an l.s.c.\ 
  mapping, and let $\psi_n\colon X\to \mathcal{P}(Y)$, $n<\omega$, be
  a sequence of mappings such that each $\psi_n$ has an
  $F_\sigma$-graph and $\psi_n(x)\cap \varphi(x)$ is a $\sigma
  Z_{n}$-set in $\varphi (x)$ for every $x\in X_n$ and $n<\omega$.
  Then, the set
  $ \{f\in \sel(\varphi): f|X_n\in
  \mathcal{M}(\psi_n|X_n)$ for every $n<\omega\}$
  is a dense
  $G_{\delta}$-set in $\sel(\varphi)$ equipped with the fine topology.
\end{theorem}

\begin{proof}
  By Theorem \ref{th:fdver}, each $\mathcal{M}(\psi_n|X_n)\cap
  \sel(\varphi|X_n)$, $n<\omega$, is a dense $G_{\delta}$-subset in
  $\sel(\varphi|X_n)$ with respect the fine topology. Let
  $\pi_n\colon\sel(\varphi)\to\sel(\varphi|X_n)$ be the map defined by
  $\pi_n(f)=f|X_n$, $f\in\sel(\varphi)$.  Note that each $\pi_n$,
  $n<\omega$, is continuous when both $\sel(\varphi)$ and
  $\sel(\varphi|X_n)$ are equipped with the fine topology. Moreover,
  $\pi_n$ is surjective because every partial selection for $\varphi
  |X_n$ can be extended to a selection for $\varphi$ (see the proof of
  \cite[Theorem 3.2$''$]{michael:56a}, also \cite{repovs-semenov:98}).
  Let us observe that each $\pi_n$, $n<\omega$, is an open map as
  well.  Towards this end, take a convex metric $\rho$ on $Y$ which is
  possible because $Y$ is a metrizable locally convex topological
  vector space, see \cite{henderson:75}. Now, it suffices to show that
  $\pi_n(V(f,\alpha)\cap\sel(\varphi))$ is open in $\sel(\varphi|X_n)$
  for every $f\in\sel(\varphi)$ and $\alpha\in C(X,(0,\infty))$.  Let
  $h\in\pi_n(V(f,\alpha)\cap\sel(\varphi))$. Then, $h=g|X_n$ for some
  $g\in V(f,\alpha)\cap\sel(\varphi)$.  Let
  $\delta(x)=\left(\alpha(x)-\rho(f(x),g(x))\right)/2$ and
  $\beta(x)=\rho(f(x),g(x))+\delta (x)$, $x\in X$. We are going to
  show that
  \[
  V(h,\delta)\cap\sel(\varphi|X_n)\subset
  \pi_n(V(f,\alpha)\cap\sel(\varphi)).
  \] 
  Indeed, for every $\ell\in V(h,\delta)\cap\sel(\varphi|X_n)$ we have
  $\rho(\ell(x),f(x))<\beta (x)$, $x\in X_n$. So, we may define an
  l.s.c.\ mapping $\phi\colon X\to\mathcal{F}_c(E)$ by $\phi
  (x)=\{\ell(x)\}$ if $x\in X_n$ and $\phi (x)=\overline{\varphi
    (x)\cap B^\rho_{\beta(x)}(f(x))}$ otherwise.  By the proof of
  \cite[Theorem 3.2$''$]{michael:56a} (see, also,
  \cite{repovs-semenov:98}), $\phi$ has a continuous selection $q$.
  Then $q\in V(f,\alpha)\cap\sel(\varphi)$ and $\ell=\pi_n(q)$.
  Therefore, $V(h,\delta)\cap\sel(\varphi|X_n)$ is a neighborhood of
  $h$ in $\pi_n(V(f,\alpha)\cap\sel(\varphi))$.\smallskip
  
  We finally accomplish the proof as follows. Since each
  $V_n=\mathcal{M}(\psi|X_n)\cap \sel(\varphi|X_n)$, $n<\omega$, is a
  dense $G_{\delta}$-set in $\sel(\varphi|X_n)$ and $\pi_n$ is
  continuous and open, $U_n=\pi_n^{-1}(V_n)$ is dense and $G_{\delta}$
  in $\sel(\varphi)$. Hence, the set $\{f\in \sel(\varphi):
  f|X_n\in \mathcal{M}(\psi|X_n),\ n<\omega\}$, being the intersection
  of all $U_n$'s, is also dense and $G_{\delta}$ in $\sel(\varphi)$.
\end{proof}

Another application of our selection theorems is the following
``strongly countable-dimensional'' analogue of Ostrand's theorem
\cite{ostrand:71}, see also \cite{engelking:95}.

\begin{theorem}\label{th:ostrand}
  For a normal space $X$ and closed subsets $X_n$, $n<\omega$, the
  following are equivalent\/$:$
\begin{itemize}
\item[\emph{(a)}] $\dim (X_n)\leq n$, $n<\omega$.\vspace*{-7pt}
\item[\emph{(b)}] For every sequence $\{\gamma_n:n<\omega\}$ of
  locally finite open covers of $X$ there is a sequence
  $\{\mu_n:n<\omega\}$ of discrete open families in $X$ such that
  for each $n<\omega$ we have\/$:$\vspace*{-7pt}
\begin{itemize}
\item[\emph{(i)}] $\mu_n$ refines $\gamma_n$,
\item[\emph{(ii)}] the union of any $n+1$ families of the sequence
  $\{\mu_{k}:k<\omega\}$ constitutes a cover of $X_n$.
\end{itemize}
\end{itemize}
\end{theorem}

\begin{proof} (a)$\ \Rightarrow\ $(b).
  We follow the proof of \cite[Theorem 1.3, implication
  $S_2\Rightarrow C$]{uspenskij:98} and \cite[Theorem 1.1, implication
  (c)$\ \Rightarrow\ $(a)]{gutev-valov:00}. Take a sequence
  $\{\gamma_n:n<\omega\}$ of locally finite open covers of $X$.  We
  are going to prove that there exists a sequence $\{\mu_n:n<\omega\}$
  of locally finite disjoint families of open sets in $X$ such that
  each $\mu_n$, $n<\omega$, refines $\gamma_n$ and the union of any
  $n+1$ of them is a cover of $X_n$.  To this end, as in \cite[Theorem
  1.1, (c)$\ \Rightarrow\ $(a)]{gutev-valov:00}, for every $n$, fix a
  metric space $(M_n,d_n)$ and a continuous map $f_n:X\to M_n$ such
  that $\{f_n^{-1}(B_2^{d_n}(z)): z\in M_n\}$ is an open cover of $X$
  refining $\gamma_n$, and $d_{n}(z,t)\le 3$ for $z,t\in M_{n}$.  Next,
  consider the disjoint union $M=\bigsqcup\{M_n:n<\omega\}$ of these
  spaces $M_n$, and let $d$ be the metric on $M$ defined as
  $d|(M_n\times M_n)=d_n,\ \hbox{and $d(z,t)=3$ provided $z\in M_i$,
    $t\in M_j$ and $i\ne j$.}$ Embed $(M,d)$ isometrically into a
  Banach space $(E,d)$, where $d$ is the metric on $E$ generated by
  the norm $\|.\|$ of $E$. Let $f=\Delta\{f_n:n<\omega\}\colon X\to
  E^{\omega}$, and let $\beta f\colon \beta X\to\beta (E^\omega)$ be
  the corresponding \v Cech-Stone extension.  Consider the space
  $H=(\beta f)^{-1}(E^\omega)$, the closure $H_n$ of each set $X_n$ in
  $H$, and the map $h\colon H\to E^\omega$ defined by $h=(\beta f)|H$.
  Then $H$ is paracompact and $\dim (H_n)=\dim (X_n)\leq n$ for every
  $n$.  Moreover, $h$ is generated by a sequence of maps $h_n\colon
  H\to E$, $n<\omega$, such that each $h_n$ extends $f_n$.  For every
  $n<\omega$, define an l.s.c.\ mapping $\varphi_n:H\to
  \mathcal{F}_c(E)$ by $\varphi_n(x)=\overline{B_1^d(h_n(x))}$, $x\in
  H$. Next, define another mapping $\varphi\colon H\to
  \mathcal{F}_c(E^\omega)$ by
  $\varphi(x)=\prod\{\varphi_n(x):n<\omega\}$, $x\in H$.  It is easily
  seen that $\varphi$ is l.s.c.\smallskip
  
  As in the proof in \cite[Theorem 1.3]{uspenskij:98} (implication
  $S_2\Rightarrow C$), we now fix a closed nowhere dense set $A\subset
  E$ such that the family $\lambda$ of all components of $E\backslash
  A$ consists of disjoint open cells of diameter $\leq 1$ and, for
  every $n<\omega$, the family $g^{-1}(\lambda)$ refines $\gamma _n$
  provided $g\colon X\to E$ is a map with $d(f_n(x),g(x))\le 1$ for
  all $x\in X$. Whenever $n<\omega$, let $\Omega_{n}=\{P\subset
  \omega: |P|=n+1\}$ and, for every $P\in\Omega_n$, let $F_{(i,P)}=A$
  if $i\in P$ and $F_{(i,P)}=E$ otherwise. Finally, set
  $F_P=\prod\{F_{(i,P)}:i<\omega\}$.  Then each $F_{P}\cap\varphi (x)$
  is, in fact, the following (possibly empty) product   
  \[
  \prod \left\{A\cap\overline{ B_1^d(f_i(x))}:i\in
    P\right\}\times\prod\left\{\overline{B_1^d(f_i(x))}:i\in
    \omega\setminus P\right\}.
  \]  
  On the other hand, each set $A\cap \overline{B_1(f_i(x))}$,
  $i<\omega$, is closed and nowhere dense in $\overline{B_1(f_i(x))}$.
  Therefore, by \cite[Corollary 2]{banakh-trushchak:99}, for every
  $n<\omega$ and $P\in\Omega_n$, the product $\prod
  \left\{A\cap\overline{B_1(f_i(x))}:i\in P\right\}$ is a $Z_{n}$-set
  in $\prod\{\varphi _i(x):i\in P\}$, $x\in H$. The last yields that
  $F_P\cap\varphi (x)$ is a $Z_{n}$-set in $\varphi (x)$ for all $x\in
  H$, $n<\omega$ and $P\in\Omega_n$. We may now apply Theorem
  \ref{th:count-case} (with $X$ replaced by $H$, $X_n$ by $H_n$, $Y$
  by $E^{\omega}$ and $\psi_n(x)=F_n=\bigcup\{F_P:P\in\Omega_n\}$,
  $x\in H$) to obtain a continuous selection $g:H\to E^{\omega}$ for
  $\varphi$ such that $g(H_n)\cap F_n=\emptyset$, $n<\omega$.  Let
  $g=\Delta\{g_n:n<\omega\}$, where each $g_n$ is a continuous map
  from $H$ into $E$. Then, $d(f_n(x),g_n(x))\le 1$ for all $x\in H$
  and $n<\omega$. Hence, according to the properties of the set $A$,
  $\mu _n=g_n^{-1}(\lambda)\cap X$ is a disjoint open family in $X$
  refining $\gamma _n$.  We can assume that each $\mu _n$ is an index
  refinement of $\gamma _n$, in particular, locally finite.  It
  remains to show that for every $n$ and $D\in \Omega_n$ the
  corresponding family $\{\mu _{i}:i\in D\}$ is a cover of $X_n$. This
  easily follows from the fact that $F_D\subset F_n$ and $g(X_n)$,
  being a subset of $g(H_n)$, avoids the set $F_n$.\medskip
  
  (b)$\ \Rightarrow\ $(a). The implication follows directly. For a fixed
  $n$, let $\alpha$ be a finite open (in $X$) cover of $X_n$. To prove
  that $\dim (X_n)\leq n$, we need to find a finite open cover
  $\lambda$ of $X_n$ refining $\alpha$ and such that any $n+2$
  elements of $\lambda$ have an empty intersection.  To this end, let
  $\gamma =\alpha\cup\{X\backslash X_n\}$. Then there exists a
  sequence $\{\mu_k:k<\omega\}$ of disjoint, open and locally finite
  families in $X$ such that each $\mu _k$ refines $\gamma$ and any
  $n+1$ of them constitute a cover of $X_n$.  We can suppose that
  every $\mu_k$ is an index refinement of $\gamma$, in particular,
  finite. Set $\lambda_{k}=\{U\in \mu_{k}:U\cap X_{n}\ne \emptyset\}$,
  $k<\omega$. Then $\lambda =\bigcup\{\lambda_k: 0\le k\le n\}$ is as
  required.
\end{proof}\medskip

By taking $X_{n}=X$, $n<\omega$, in Theorem \ref{th:ostrand} we now
get the following ``finite-dimensional'' generalization of the
Ostrand's theorem (recall that Ostrand originally proved his theorem
for finite-dimensional spaces).

\begin{corollary}\label{cor:ostrand}
  For a normal space $X$ and $n<\omega$ the  following are equivalent\/$:$
\begin{itemize}
\item[\emph{(a)}] $\dim (X)\leq n$.\vspace*{-7pt}
\item[\emph{(b)}] For every sequence $\{\gamma_k:k<\omega\}$ of
  locally finite open covers of $X$ there is a sequence
  $\{\mu_k:k<\omega\}$ of discrete open families in $X$ such that each
  $\mu_k$, $k<\omega$, refines $\gamma_k$ and the union of any $n+1$
  elements of the sequence $\{\mu_k:k<\omega\}$ constitutes a cover of
  $X$.
\end{itemize}
\end{corollary}

\section{Closed maps and selections}

Let $M$ be a metrizable space. We will use the term \emph{distance}
for a possibly infinite-valued function $d:M\times M\to [0,+\infty]$
which satisfies the axioms of a metric on $M$. So, let
$\mathcal{D}_{\infty}(M)$ denote the set of all distances on $M$ which
are compatible with the topology of $M$. Then, to every $d\in
\mathcal{D}_{\infty}(M)$ and every positive real number $r\in\R$ we
may associate a metric $d_{r}\in \mathcal{D}(M)$ on $M$ defined by
$d_{r}(x,y)=\min\{r, d(x,y)\}$, $x,y\in M$. Now, we shall say that
$d\in\mathcal{D}_\infty(M)$ has a property $\mathcal{P}$ if $d_r$ has
$\mathcal{P}$ for every $r>0$. For instance, \emph{completeness} is
such a property $\mathcal{P}$.\medskip

Here is a natural example of distances we will deal with.
Let $X$ be a space, and let $(Y,d)$ be a metric space. In what
follows, we will rely on the \emph{uniform topology} on $C(X,Y)$
generated by the distance
\[
\bfd(f,g)=\sup\{d(f(x),g(x)): x\in X\},\qquad f,g\in C(X,Y).
\]
Note that $\bfd_r(f,g)=\sup\{d_r(f(x),g(x)): x\in X\}$, whenever
$r>0$. Thus, $(C(X,Y),\bfd)$ is complete if $d$ is a complete metric
on $Y$ (see, for instance, \cite[Chapter XII]{dugundji:66}).
Concerning maps $\ell:Z\to C(X,Y)$, let us agree to denote by
$\ell[z]\in C(X,Y)$ the value of $\ell$ in a particular point $z\in
Z$.  \medskip

Throughout this section, to every surjective map  $f:X\to Y$ and a
set-valued mapping  $\varphi:X\to 2^E$ we associate another set-valued
mapping $\Delta_{(f,\varphi)}:Y\to \mathcal{P}(C(X,E))$ defined by 
\[
\Delta_{(f,\varphi)}(y)=\{g\in C(X,E):
g|f^{-1}(y)\in\sel(\varphi|f^{-1}(y))\},\quad y\in Y.
\] 
Also, let us recall that a subset $A$ of a space $X$ is
\emph{$P$-embedded} in $X$ if every continuous map $g:A\to E$ into a
Banach space $E$ is continuously extendable over $X$, see
\cite{alo-sennott:71,morita:75b,przymusinski:78,shapiro:66}. Finally,
by a Banach space $(E,d)$ we mean a Banach space $E$ and a metric $d$
generated by the norm of $E$.

\begin{theorem}\label{th:constr}
  Let $X$ and $Y$ be a spaces, and let $f:X\to Y$ be a continuous
  closed surjection such that each $f^{-1}(y)$, $y\in Y$, is
  paracompact and $P$-embedded in $X$. Also, let $(E,d)$ be a Banach
  space and $\varphi:X\to \mathcal{F}_c(E)$ be an l.s.c.\ mapping.
  Then, $\Delta_{(f,\varphi)}$ is an l.s.c.\ mapping from $Y$ with
  values in $\mathcal{F}_c(C(X,E))$ provided $C(X,E)$ is endowed with
  the uniform topology generated by the distance $\bfd$.  Moreover, if
  $\ell$ is a continuous selection for $\Delta_{(f,\varphi)}$, then
  $g(x)=\ell[f(x)](x)$, $x\in X$, defines a continuous selection
  for~$\varphi$.
\end{theorem}

\begin{proof}
  Note that, by \cite[Theorem 3.2$''$]{michael:56a},
  $\sel(\varphi|f^{-1}(y))\ne\emptyset$ for every $y\in Y$ because
  each $f^{-1}(y)$ is paracompact. Hence,
  $\Delta_{(f,\varphi)}(y)\ne\emptyset$ for every $y\in Y$ because
  each $f^{-1}(y)$ is $P$-embedded in $X$. Thus,
  $\Delta_{(f,\varphi)}:Y\to \mathcal{F}_c(C(X,E))$ because
  $\varphi:X\to \mathcal{F}_c(E)$.\smallskip
  
  To show that $\Delta_{(f,\varphi)}$ is l.s.c. take a point $y_0\in
  \Delta_{(f,\varphi)}^{-1}(B^{\bfd}_{\varepsilon}(g_0))$, where
  $g_0\in\Delta_{(f,\varphi)}(y_0)$ and $\varepsilon>0$. It suffices
  to find a neighborhood of $y_0$ in $Y$ which is contained in
  $\Delta_{(f,\varphi)}^{-1}(B^{\bfd}_{\varepsilon}(g_0))$. Since
  \[
  V=\{x\in X:
  d(g_0(x),\varphi(x))<\varepsilon/2\}
  \]
  is a neighbourhood of $f^{-1}(y_0)$ and $f$ is closed, there exists
  of a neighbourhood $U$ of $y_0$, with $f^{-1}(U)\subset V$.  Let
  show that
  $U\subset\Delta_{(f,\varphi)}^{-1}(B^{\bfd}_{\varepsilon}(g_0))$.
  Indeed, take an $y\in U$ and define an l.s.c.\ mapping $\Phi:
  f^{-1}(y)\to \mathcal{F}_c(E)$ by
  $\Phi(x)=\overline{B^d_{\varepsilon/2}(g_0(x))\cap \varphi(x)}$,
  $x\in f^{-1}(y)$.  Then, by \cite[Theorem 3.2$''$]{michael:56a},
  $\Phi$ has a continuous selection $h$. Thus, $\varphi|f^{-1}(y)$ has
  a continuous selection $h$, with
  $\bfd(h,g_0|f^{-1}(y))\le\varepsilon/2<\varepsilon$.

  \begin{claim}
    There exists a continuous extension $g:X\to E$ of $h$ such that
    $\bfd(g,g_0)<\varepsilon$. 
  \end{claim}
  
  In order to show our claim, we first extend $h$ to a continuous map
  $g_{1}:X\to E$ which is possible because $f^{-1}(y)$ is $P$-embedded
  in $X$. Next, let us observe that $Z=\{x\in X:
  d(g_{0}(x),g_{1}(x))\ge 2\varepsilon/3\}$ is a zero-set of $X$, with
  $Z\cap f^{-1}(y)=\emptyset$.  Since $f^{-1}(y)$ is $P$-embedded, by
  \cite[Corollary 3.6.B]{blair-hager:74}, there now exists a
  continuous function $\eta:X\to \uin$ such that $f^{-1}(y)\subset
  \eta^{-1}(0)$ and $Z\subset \eta^{-1}(1)$.  Finally, we may define
  $g:X\to E$ by $g(x)=\eta(x)\cdot g_{0}(x)+(1-\eta(x))\cdot
  g_{1}(x)$, $x\in X$.  This $g$ is as required.\medskip
 
  Now, we have that $g\in \Delta_{(f,\varphi)}(y)$ because $g$ is an
  extension of $h$, while, by Claim~1, $y\in
  \Delta_{(f,\varphi)}^{-1}(B^{\bfd}_{\varepsilon}(g_0))$, so
  $U\subset \Delta_{(f,\varphi)}^{-1}(B^{\bfd}_{\varepsilon}(g_0))$.
  Therefore, $\Delta_{(f,\varphi)}$ is l.s.c.\medskip
  
  Let prove the final part of Theorem 6.1. Suppose $\ell$ is a
  continuous selection for $\Delta_{(f,\varphi)}$, and let
  $g(x)=\ell[x](x)$ for every $x\in X$. Clearly $g$ is a selection for
  $\varphi$. So, it only remains to show that $g$ is continuous.  Take
  a point $x_0\in X$ and $\varepsilon>0$, and let
  $U=\ell^{-1}(B_{\varepsilon/2}^{\bfd}(\ell[x_0]))\cap
  (\ell[x_0])^{-1}(B_{\varepsilon/2}^d(\ell[x_0](x_0)))$.  Then, $x\in
  U$ implies
  \begin{eqnarray*}
    d(g(x),g(x_0)) &=& d(\ell[x](x),\ell[x_0](x_0))\\
    &\le&  d(\ell[x](x),\ell[x_0](x)) +
    d(\ell[x_0](x),\ell[x_0](x_0))\\
    &\le&  \bfd(\ell[x],\ell[x_0]) +
    d(\ell[x_0](x),\ell[x_0](x_0))\\
    &<& \varepsilon/2 +\varepsilon/2=\varepsilon.
  \end{eqnarray*}
  This completes the proof.
\end{proof}

Theorem \ref{th:constr} may have some general interest. For instance,
it implies the following slight generalization of a result of Hanai
\cite{hanai:56}, see also \cite{engelking:89}.

\begin{theorem}\label{th:hanai}
  Let $X$ be a space, $Y$ be a paracompact space, and $f:X\to Y$ be a
  continuous closed surjection such that each $f^{-1}(y)$, $y\in Y$,
  is paracompact and $P$-embedded in $X$. Then, $X$ is paracompact
  too. 
\end{theorem}

\begin{proof}
  Take an l.s.c.\ mapping $\varphi:X\to \mathcal{F}_c(E)$ for some
  Banach space $(E,d)$. By \cite[Theorem 3.2$''$]{michael:56a}, it
  suffices to show that $\varphi$ has a continuous selection. Consider
  the l.s.c.\ mapping $\Delta_{(f,\varphi)}:Y\to
  \mathcal{F}_c(C(X,E))$, where $C(X,E)$ is endowed with the uniform
  topology generated by the distance $\bfd$. In fact, the same
  topology on $C(X,E)$ is generated by the metric $\bfd_1$ on $C(X,E)$
  defined by $\bfd_1(g,h)=\min\{1,\bfd(g,h)\}$, see the beginning of
  the section. Since $(C(X,E),\bfd_1)$ is a Fr{\'e}chet space,
  $\Delta_{(f,\varphi)}$ has a continuous selection \cite{michael:56a}
  (see, also, \cite{repovs-semenov:98}). Hence, Theorem
  \ref{th:constr} completes the proof.
\end{proof}

Finally, we provide an alternative proof of the dimension-lowering
mapping theorem, see \cite{engelking:95} for the proof and history of
this theorem.

\begin{theorem}\label{th;hurewicz}
  If $f\colon X\to Y$ is a closed continuous mapping of a normal space
  $X$ to a paracompact space $Y$ and there exists an integer $k\geq 0$
  such that $\dim(f^{-1}(y))\leq k$ for every $y\in Y$, then
  $\dim(X)\leq\dim(Y)+k$.
\end{theorem}

\begin{proof}
  We may assume that $X$ is paracompact. Indeed, we may consider the
  {\v C}ech-Stone extension $\beta f: \beta X\to \beta Y$. Then, by
  Theorem \ref{th:hanai}, $H=(\beta f)^{-1}(Y)$ is a paracompact
  space, while $\tilde{f}=\beta f|H$ is a perfect map such that
  $\tilde{f}^{-1}(y)=\overline{f^{-1}(y)}$ for every $y\in Y$. Since
  $X$ is normal, we have that $\beta(f^{-1}(y))=\tilde{f}^{-1}(y)$, so
  $\dim(\tilde{f}^{-1}(y))\le k$ for every $y\in Y$. Finally, $\beta
  H=\beta X$ implies that $\dim(H)=\dim(X)$.\smallskip
  
  Thus, let $X$ and $Y$ be paracompact, and let $f:X\to Y$ be a closed
  continuous surjection as in the hypothesis.  Also, suppose that
  $\dim(Y)\leq m$.  To show that $\dim (X)\leq m+k$, it suffices to
  show that every map $g\colon X\to\R^{m+k+1}$ is removable from the
  origin $\vartheta$ of $\R^{m+k+1}$. For a fixed $g_0\in
  C(X,\R^{m+k+1})$ and $\varepsilon >0$, we define a set-valued mapping
  $\varphi:X\to \mathcal{F}_c(\R^{m+k+1})$ by
  $\varphi(x)=\overline{B^d_{\varepsilon}(g_0(x))}$, $x\in X$, where $d$
  is the usual Euclidean metric on $\R^{m+k+1}$. Thus, by Theorem
  \ref{th:constr}, we may get an l.s.c.\ mapping
  $\Phi=\Delta_{(f,\varphi)}:Y\to
  \mathcal{F}_c(C(X,\R^{m+k+1}))$.\smallskip
  
  Now we consider the mapping $\Psi\colon Y\to
  \mathcal{F}(C(X,\R^{m+k+1}))$ defined by
  \[
  \Psi (y)=\{g\in C(X,\R^{m+k+1}):\vartheta\in g(f^{-1}(y))\}, \quad
  y\in Y.
  \]
  It is a routine verification that the graph of $\Psi$ is closed. Let
  us show that $\Psi (y)\cap\Phi(y)$ is a $Z_m$-set in $\Phi(y)$ for
  every $y\in Y$. Take a fixed point $y\in Y$ and a continuous map
  $u\colon\ball^m\to\Phi (y)$. We are going to prove that for every
  $\delta >0$ there exists a continuous map $v\colon\ball ^m\to\Phi
  (y)$ which is $\delta$-close to $u$ and avoids the set $\Psi (y)$.
  Observe that $u$ generates a continuous map $p\colon\ball ^m\times
  X\to\R^{m+k+1}$, $p(z,x)=u[z](x)$, such that $p(z,x)\in\varphi (x)$
  for every $x\in f^{-1}(y)$ and $z\in\ball^m$. Define
  $\phi\colon\ball ^m\times f^{-1}(y)\to \mathcal{F}_c(\R^{m+k+1})$ by
  $\phi(z,x)=\overline{B^d_{\delta/2}(p(z,x))\cap\varphi (x)}$,
  $(z,x)\in\ball^m\times f^{-1}(y)$. Since $\vartheta$ is a
  $Z_{m+k}$-set in each $\varphi (x)$ (as a $Z_{m+k}$-set in
  $\R^{m+k+1}$), it is a $Z_{m+k}$-set in $\phi(z,x)$ for any
  $(z,x)\in\ball^m\times f^{-1}(y)$ (see \cite[Lemma
  2.3]{gutev-valov:00}). Moreover, $\phi$ is an l.s.c.\ and
  $\dim(\ball ^m\times f^{-1}(y))\leq m+k$. Hence, by Theorem
  \ref{th:fdver}, $\phi$ has a continuous selection
  $\tilde{q}\colon\ball^m\times f^{-1}(y)\to\R^{m+k+1}$ avoiding
  $\vartheta$. In particular, we get that $\tilde{q}(z,x)\in
  \varphi(x)$ for every $x\in f^{-1}(y)$ and $z\in \ball^m$, while
  $\bfd(\tilde{q},p|\ball^m\times f^{-1}(y))\le\delta/2$.  So, we may
  extend $\tilde{q}$ to a continuous map $q:\ball^m\times
  X\to\R^{m+k+1}$ such that $\bfd(q,p)<\delta$. Then $q$ determines a
  continuous map $v\colon\uin ^m\to C(X,\R^{m+k+1})$, with
  $v[z](x)=q(z,x)$.  This $v$ is as required.\smallskip
  
  We now finish the proof as follows. By Theorem \ref{th:fdver},
  $\Phi$ has a continuous selection $\ell\colon Y\to C(X,\R^{m+k+1})$
  such that $\ell(y)\in\Phi (y)\backslash\Psi (y)$ for every $y\in Y$.
  Then the map $g\colon X\to\R^{m+k+1}$, $g(x)=\ell[f(x)](x)$, is
  $\varepsilon$-close to $g_0$ and avoids $\vartheta$.
\end{proof}


\newcommand{\noopsort}[1]{} \newcommand{\singleletter}[1]{#1}
\providecommand{\bysame}{\leavevmode\hbox to3em{\hrulefill}\thinspace}
\providecommand{\MR}{\relax\ifhmode\unskip\space\fi MR }
\providecommand{\MRhref}[2]{%
  \href{http://www.ams.org/mathscinet-getitem?mr=#1}{#2}
}
\providecommand{\href}[2]{#2}

\bigskip\bigskip

{\footnotesize\sc School of Mathematical and Statistical Sciences,
  Faculty of Science, University of Natal, King George V Avenue,
  Durban 4041, South Africa}

{\em E-mail address\/}: {\tt gutev@nu.ac.za}\bigskip

{\footnotesize\sc Department of Computer Science and Mathematics, Nipissing University, 100
  College Drive, P.O. Box 5002, North Bay, ON, P1B 8L7, Canada}

{\em E-mail address\/}: {\tt veskov@unipissing.ca}
\end{document}